\documentclass[11pt,fleqn,twoside]{article}
\usepackage{amsfonts,amssymb,latexsym}
\makeatletter
\newcommand{\prava}[1]{\small\it
\begin{flushleft}
Copyright \copyright \ 1999 by  #1
\end{flushleft}}

\newcommand{\name}[1]{\begin{flushleft}
                       \LARGE \bf #1
                       \end{flushleft}\vspace{-3mm}}

\newcommand{\Author}[1]{\begin{flushleft}
                       \it #1 \end{flushleft}}

\newcommand{\Adress}[1]{\begin{flushleft}
                       \it #1 \end{flushleft}}

\newcommand{\Date}[1]{\begin{flushleft}
                      \small  \it #1 \end{flushleft}}

\newcommand{\ehkol}{Author \ name}
\newcommand{\ohkol}{Article \ name}
\renewcommand{\@evenhead}{
\hspace*{-3pt}\raisebox{-15pt}[\headheight][0pt]{\vbox{\hbox to \textwidth 
{\thepage \hfil \ehkol}\vskip4pt \hrule}}}
\renewcommand{\@oddhead}{
\hspace*{-3pt}\raisebox{-15pt}[\headheight][0pt]{\vbox{\hbox to \textwidth 
{\ohkol \hfil \thepage}\vskip4pt\hrule}}}
\renewcommand{\@evenfoot}{}
\renewcommand{\@oddfoot}{}

\newcommand{\be}{\begin{equation}}
\newcommand{\ee}{\end{equation}}
\newcommand{\ba}{\hspace*{-5pt}\begin{array}}
\newcommand{\ea}{\end{array}}

\newcommand{\ds}{\displaystyle}
\makeatother

\usepackage{lscape}

\begin{document}
\thispagestyle{empty}
\setcounter{page}{294}
\renewcommand{\ehkol}{P. Bracken and A.M. Grundland}
\renewcommand{\ohkol}{On Certain Classes of Solutions of the Weierstrass-Enneper
System} 

\begin{flushleft}
\footnotesize \sf
Journal of Nonlinear Mathematical Physics \qquad 1999, V.6, N~3,
\pageref{bracken-fp}--\pageref{bracken-lp}.
\hfill {\sc Article}
\end{flushleft}

\vspace{-5mm}

\renewcommand{\footnoterule}{}
{\renewcommand{\thefootnote}{} 
 \footnote{\prava{P. Bracken and A.M. Grundland}}}

\name{On Certain Classes of Solutions of \\
the Weierstrass-Enneper System Inducing \\
Constant Mean Curvature Surfaces}\label{bracken-fp}

\Author{P. BRACKEN~$^{\dag \ddag}$ and A. M. GRUNDLAND~$^{\dag}$}

\Adress{$\dag$~Centre de Recherches Math\'{e}matiques, Universit\'{e} de Montr\'{e}al, \\
~~2920 Chemin de la Tour, Pavillon Andr\'{e} Aisenstatt, C. P. 6128 Succ. Centre Ville, \\
~~Montr\'{e}al, QC, H3C 3J7 Canada\\[1mm]
$\ddag$~Department of Mathematics and Statistics, McGill University,  \\
~~Montr\'{e}al, QC, H3A 2K6, Canada}

\Date{Received December 28, 1998; Accepted March 15, 1999}

\begin{abstract}
\noindent
Analysis of the generalized Weierstrass-Enneper system includes
the estimation of the degree of indeterminancy of the general
analytic solution and the discussion of the boundary value
problem. Several dif\/ferent procedures for constructing
certain classes of solutions to this system, including potential,
harmonic and separable types of solutions, are proposed.
A technique for reduction of the Weierstrass-Enneper
system to decoupled linear equations, by subjecting it to
certain dif\/ferential constraints, is presented as well.
New elementary and doubly periodic solutions are found,
among them kinks, bumps and multi-soliton solutions.
\end{abstract}

\renewcommand{\theequation}{\arabic{section}.\arabic{equation}}


\section{Introduction}

The Weierstrass-Enneper system [1] has proved to be a very useful and
suitable tool in the study of minimal surfaces in ${\mathbb R}^{3}$.
Since Weierstrass and Enneper, this subject has been investigated
extensively by many authors (eg. Kenmotsu~[2],
Hof\/fman and Osserman~[3], Konopelchenko~[4--9]).

The original formulation by Weierstrass and Enneper [1]
of a system inducing minimal surfaces
can be presented brief\/ly as follows. Let $\alpha$ and $\beta$ be 
holomorphic functions that satisfy $\bar{\partial} \alpha =0$,
and $\bar{\partial} \beta = 0$ such that the equations
\[
\partial w_{1} = i (\alpha^{2} + \beta^{2} ), \qquad
\partial w_{2} = \alpha^{2} - \beta^{2},   \qquad
\partial w_{3} = - 2 \alpha \beta,
\]
hold, where the derivatives are abbreviated $\partial = \partial / \partial
z$, $\bar{\partial} = \partial / \partial \bar{z}$. Introduce a system of
three real-valued functions $X_{i} (z, \bar{z} )$, $i=1,2,3$ which can be
considered a coordinate system for a surface embedded in
${\mathbb R}^{3}$, def\/ined as follows
\be
\ba{l}
\ds X_{1} = Re \, w_{1} = Re \, \left[ i \int_{C} (\alpha^{2} + \beta^{2} ) \, dz' \right],
\vspace{3mm}\\
\ds X_{2} = Re \, w_{2} = Re \, \left[ \int_{C} (\alpha^{2} - \beta^{2} ) \, dz' \right],
\vspace{3mm}\\
\ds X_{3} = Re \, w_{3} = - Re \, \left[ 2 \int_{C} \alpha \beta \, dz' \right],
\ea
\ee
where $C$ is any contour in the common domain of holomorphicity of the
functions $\alpha$ and~$\beta$.
The $X_{i}$ def\/ine a minimal surface with $z=c_{1}$ and $\bar{z} = c_{2}$
as minimal coordinate lines on this surface, respectively.

More recently, the generalized Weierstrass-Enneper (WE)
representation for constant mean curvature surfaces in ${\mathbb R}^{3}$
has been introduced by B.~Konopelchenko [7--9] and his formula
is a starting point for our analysis.
Namely, we consider the nonlinear system, a
type of two-dimensional Dirac equation
for the f\/ields $\psi_{1}$ and $\psi_{2}$, given by
\be
\partial \psi_{1} = p \psi_{2}, \qquad  \bar{\partial} \psi_{2} =- p \psi_{1},  \qquad
\bar{\partial} \bar{\psi}_{1} = p \bar{\psi}_{2}, \qquad 
\partial \bar{\psi}_{2} = - p \bar{\psi}_{1},
\ee
where
\[
p = |\psi_{1}|^{2} + |\psi_{2}|^{2},
\]
in the neighbourhood of some point $(z_{0}, \bar{z}_{0}) \in {\mathbb C}$,
and the bar on $\psi_{i}$ denotes complex conjugation. 
In this paper, when we refer to the WE system, we mean the
modif\/ied version (1.2) of the original Weierstrass Enneper system.

One then def\/ines the three real-valued functions $X_{i} (z, \bar{z})$, 
$i=1,2,3$, by means of the formulae
\be
\ba{l}
\ds X_{1} + i X_{2} = 2i \int_{\gamma} \left(\bar{\psi}_{1}^{2} \, dz' -
\bar{\psi}_{2}^{2} \, d \bar{z}'\right),
\vspace{3mm}\\
\ds X_{1} - i X_{2} = 2i \int_{\gamma} \left(\psi_{2}^{2} \, dz' - 
\psi_{1}^{2} \, d \bar{z}'\right),
\vspace{3mm}\\
\ds X_{3} = -2 \int_{\gamma} \left( \bar{\psi}_{1} \psi_{2} \, dz' +
\psi_{1} \bar{\psi}_{2} \, d \bar{z}'\right),
\ea
\ee
where $\gamma$ is any contour in ${\mathbb C}$. On account of the
system (1.2), the right hand side of~(1.3) does not depend on the
choice of $\gamma$.
It was shown in~[4,~5] that for each pair of solutions
$( \psi_{1}, \psi_{2})$ of WE system (1.2), the 
formulae (1.3) determine a conformal immersion
of a constant mean curvature surface in ${\mathbb R}^{3}$.
The induced metric on the surface and its
Gaussian curvature are given by~[7,~14]
\be
ds^{2} = 4 p^{2} \, dz \, d \bar{z}, \qquad
K = - p^{-2} \partial \bar{\partial} (\ln p)
\ee 
in isothermic coordinates, respectively.
Finally, a well known property of WE system (1.2), in the context
of the sigma model~[10--12], is the existence of a topological charge
\be
I = - \frac{i}{2 \pi} \int_{\gamma} \frac{1}{p^{2}} \left(|j|^{2} - p^{4}\right) \,
dz d \bar{z},
\ee
where $j$ is an entire function[14] def\/ined by 
\be
j(z) = \bar{\psi}_{1} \partial \psi_{2} - \psi_{2} \partial \bar{\psi}_{1}.
\ee
In fact, $j$ is a conserved quantity since
\be
\bar{\partial} j = \bar{\partial} \left(\bar{\psi}_{1} \partial \psi_{2} -
\psi_{2} \partial \bar{\psi}_{1} \right) = 0
\ee
holds, whenever the WE system (1.2) is satisf\/ied. The scalar function
$j(z)$ is referred to as the current for WE system (1.2)~[7].
Note that if the integral in (1.5) exists, then $I$ is an integer.

In this paper we explore several dif\/ferent aspects of the
generalized WE system (1.2). We investigate certain general
characteristics of this system and propose several new 
approaches to the construction of its solutions.
The paper is organized as follows. Section~2
contains a detailed account of the estimation of the degree of 
freedom of the general analytic solution to the WE system,
based on Cartan's theory of systems in involution.
In Section~3 we discuss the existence and
uniqueness of solutions to
a boundary value problem for the WE system. 
Section~4 and Section~5 contain several examples of potential and
harmonic solutions of the WE system which include an
explicit form of an algebraic multi-soliton solution. 
Next, in Section~6 we introduce a set of 
dif\/ferential constraints under which the WE system can be
reduced to a system of linear coupled equations
and we construct several examples of solutions using this approach. 
Section~7 presents a certain variant of the separation of variables 
technique applied to the WE system 
which allows us to construct solitonlike solutions (bumps and kinks).

\setcounter{equation}{0}

\section{The estimation of the degree of indeterninancy of the ge\-ne\-ral
analytic solution to the Weierstrass-Enneper system}

Now, using Cartan's theorem on systems in involution [13], we estimate the
degree of indeterminancy of the general analytic solution of WE system (1.2).
For this purpose, we perform the analysis using the apparatus of 
dif\/ferential forms which are equivalent to the initial system.
The problem is reduced to examining the Cartan numbers of these
exterior forms and the number of arbitrary parameters admitted 
by the general solution of the system of polar equations.

For computational purposes, it is useful to introduce the following notation
\be
\ba{l}
x = (x^{1},x^{2}):=(\bar{z},z), \qquad
\xi = (\xi^{1},\xi^{2},\xi^{3},\xi^{4}):=
(\bar{\partial} \psi_{1}, \partial \psi_{2}, \partial \bar{\psi}_{1},
\bar{\partial} \bar{\psi}_{2}),
\vspace{2mm}\\
u= (u_{1},u_{2},u_{3},u_{4},u_{5},u_{6},u_{7},u_{8}):=
(\psi_{1}, \psi_{2}, \bar{\psi}_{1},\bar{\psi}_{2}, \partial \psi_{1},
\bar{\partial} \psi_{2}, \bar{\partial} \bar{\psi}_{1},
\partial \bar{\psi}_{2}) .
\ea
\ee
It means that we interpret $z$ and $\bar{z}$ as independent 
coordinates $x^{1}$ and $x^{2}$, respectively, in ${\mathbb R}^{2}$ space, 
and the coordinates $(u_{1}, \ldots , u_{8})$
as independent variables in ${\mathbb R}^{8}$ space.
The quantity $\xi$ represents a vector of all f\/irst derivatives 
of $\psi_{i}$ which do not
appear in the WE system. Under notation (2.1)
the system (1.2) becomes
\be
u_{5} = p u_{2}, \qquad u_{7} = p u_{4}, \qquad
u_{6} = -p u_{1}, \qquad  u_{8} = - p u_{3}, \qquad
p = u_{1}u_{3} + u_{2} u_{4},
\ee
and the dif\/ferentiation of $p$
with respect to $z$ and $\bar{z}$ yields
\be
\partial p = u_{1} \xi^{3} + u_{4} \xi^{2}, \qquad
\bar{\partial} p = u_{3} \xi^{1} + u_{2} \xi^{4},
\ee
whenever (1.2) holds.
Note that the $\xi^{s}$ enter linearly into the expressions $\partial p$
and $\bar{\partial} p$. If we consider the variables
$u= (u_{1}, \ldots , u_{8})$ and $\xi = (\xi^{1}, \ldots , \xi^{4})$
as unknown functions of $x = (x^{1}, x^{2})$ then, in terms of
(2.1) and (2.2), the WE system (1.2) is equivalent to the system
of dif\/ferential one-forms
\be
\ba{l}
\omega_{1} = du_{1} - (\xi^{1} \, dx^{1} + p u_{2} \, dx^{2}) =0,
\vspace{2mm}\\
\omega_{2} = du_{2} - (-p u_{1} \, dx^{1} + \xi^{2} \, dx^{2}) =0,
\vspace{2mm}\\
\omega_{3} = du_{3} - (p u_{4} \, dx^{1} + \xi^{3} \, dx^{2}) =0,
\vspace{2mm}\\
\omega_{4} = du_{4} - (\xi^{4} \, dx^{1} -p u_{3} \, dx^{2}) = 0,
\vspace{2mm}\\
\omega_{5} = du_{5} -[ u_{2} \bar{\partial} p - p^{2} u_{1}] \, dx^{1}
- [ u_{2} \partial p + p \xi^{2} ] \, dx^{2} =0,
\vspace{2mm}\\
\omega_{6} = du_{6} - [ u_{1} \bar{\partial} p + p \xi^{1} ] \, dx^{1}
- [ u_{1} \partial p + p^{2} u_{2} ] \, dx^{2} =0,
\vspace{2mm}\\
\omega_{7} = du_{7} - [ u_{4} \bar{\partial} p + p \xi^{4}] \, dx^{1}
- [ u_{4} \partial p - p^{2} u_{3} ] \, dx^{2} = 0,
\vspace{2mm}\\
\omega_{8} = du_{8} - [ u_{3} \bar{\partial} p + p^{2} u_{4} ] \, dx^{1}
- [ u_{3} \partial p + p \xi^{3} ] \, dx^{2} = 0,
\ea
\ee
in two independent variables $x^{1}$, $x^{2}$ which form some local coordinate
system in the real space ${\mathbb R}^{2}$.
System (2.4) can be written in the abbreviated form
\be
\omega_{s} = du_{s} - G_{s \mu} (x, \xi, u) \, dx^{\mu}, \qquad
s=1, \ldots , 8, \quad \mu=1,2,
\ee
where the functions $G_{s \mu}$ depend only on $(x,\xi, u)$
and where $\xi$ enters linearly into $G_{s \mu}$, due to (2.3).
We are interested in the evaluation of the degree of freedom
of the most general analytic solution of (2.4)
which can be expressed by
\[
u_{s} = u_{s} (x^{1}, x^{2}), \qquad \xi^{r} = \xi^{r} (x^{1}, x^{2}),
\qquad s=1, \ldots, 8, \quad r=1, \ldots , 4.
\]
According to Cartan's Theorem on systems in involution [13], we can
formulate the following. 

\medskip

\noindent
{\bf Proposition 1.} {\it If the system of
differential one-forms (2.4) is in involution at a regular
point $(x_{0}, \xi_{0}, u_{0})$ and if it is an analytic system in some
neighbourhood of $(x_{0}, \xi_{0}, u_{0})$, then
the general analytic solution of (2.4) with independent variables $x^{1}$,  $x^{2}$
exists in some neighbourhood of a regular point $(x_{0}, \xi_{0}, u_{0})$
and it depends on four arbitrary real analytic functions of one 
real variable. }

\medskip

\noindent
{\bf Proof.} Under notation (2.1) and relations (2.3), the exterior dif\/ferentiation
of system~(2.4) leads to the following $2$-form system whenever
system (2.4) holds
\[
\ba{l}
\ds \Omega_{1} \equiv d \omega_{1} = dx^{1} \wedge d \xi^{1}
-\left[ u_{2} \bar{\partial}p - p^{2} u_{1}\right] \, dx^{1} \wedge dx^{2},
\vspace{3mm}\\
\ds \Omega_{2} \equiv d \omega_{2} = -\left[ u_{1} \partial p + p^{2} u_{2}\right] \,
dx^{1} \wedge dx^{2} + dx^{2} \wedge d \xi^{2},
\vspace{3mm}\\
\ds \Omega_{3} \equiv d\omega_{3} = \left[ u_{4} \partial p - p^{2} u_{3} \right] \,
dx^{1} \wedge dx^{2} + dx^{2} \wedge d \xi^{3},
\vspace{3mm}\\
\ds \Omega_{4} \equiv d \omega_{4} = dx^{1} \wedge d \xi^{4} +
\left[ u_{3} \bar{\partial} p + p^{2} u_{4} \right] \, dx^{1} \wedge dx^{2},
\vspace{3mm}\\
\ds \Omega_{5} \equiv d \omega_{5} =- u_{2} u_{3} \,
d \xi^{1} \wedge dx^{1} + u_{2}^{2} dx^{1} \wedge d \xi^{4} - u_{1} u_{2}
d \xi^{3} \wedge dx^{2} - (p + u_{2} u_{4})\, d \xi^{2} \wedge dx^{2} 
\vspace{3mm}\\
\qquad + [ u_{1}u_{2}( u_{4} \partial p - p^{2}u_{3}) -(p+u_{2}u_{4})
(u_{1} \partial p  + p^{2} u_{2})  
\vspace{3mm}\\
\ds \qquad - u_{2}u_{3} ( u_{2} \bar{\partial} p - p^{2} u_{1}) + u_{2}^{2} 
((u_{3} \bar{\partial} p + p^{2} u_{4})] \, dx^{1} \wedge dx^{2},
\ea
\]
\be
\ba{l}
\ds \Omega_{6} \equiv d \omega_{6} = -(p+u_{1}u_{3}) \, dx^{1} \wedge d \xi^{1}
-u_{1}u_{2} \, dx^{1} \wedge d \xi^{4} - u_{1} u_{4} \, dx^{2} \wedge d \xi^{2}
\vspace{3mm}\\
\ds \qquad -u_{1}^{2} \, dx^{2} \wedge d \xi^{3}
+[ - u_{1}^{2} ( u_{4} \partial p - p^{2} u_{3}) + u_{1} u_{4}
(u_{1} \partial p + p^{2} u_{2})
\vspace{3mm}\\
\ds \qquad - u_{1} u_{2} ( u_{3} \bar{\partial} p + p^{2} u_{4}) + (p+u_{1} u_{3})
(u_{2} \bar{\partial} p - p^{2} u_{1})] \, dx^{1} \wedge dx^{2},
\vspace{3mm}\\
\ds \Omega_{7} \equiv d \omega_{7} = u_{3} u_{4} \, dx^{1} \wedge d \xi^{1}+
(p+ u_{2} u_{4}) \, dx^{1} \wedge d \xi^{4} + u_{4}^{2} \, dx^{2} \wedge
d \xi^{2} 
\vspace{3mm}\\
\ds \qquad + u_{1} u_{4} \, dx^{2} \wedge d \xi^{3}
+[ u_{1}u_{4}( u_{4} \partial p - p^{2} u_{3}) -u_{4}^{2}(u_{1} \partial p
+p^{2} u_{2})
\vspace{3mm}\\
\qquad - u_{3} u_{4} ( u_{2} \bar{\partial} p - p^{2} u_{1}) +(p+ u_{2} u_{4})
( u_{3} \bar{\partial} p + p^{2} u_{4})] \, dx^{1} \wedge dx^{2},
\vspace{3mm}\\
\ds \Omega_{8} \equiv d \omega_{8} =- u_{3}^{2} \, dx^{1} \wedge  d \xi^{1}
-u_{2} u_{3} \, dx^{1} \wedge d \xi^{4} - u_{3} u_{4} \, dx^{2} \wedge
d \xi^{2} 
\vspace{3mm}\\
\ds \qquad -(p+u_{1}u_{3}) \, dx^{2} \wedge d \xi^{3}
+ [ u_{3} u_{4}( u_{1} \partial p + p^{2} u_{2}) - (p +u_{1} u_{3})
( u_{4} \partial p - p^{2} u_{3})
\vspace{3mm}\\
\ds \qquad +u_{3}^{2} ( u_{2} \bar{\partial} p - p^{2} u_{1}) - u_{2} u_{3}
(u_{3} \bar{\partial} p + p^{2} u_{4} )] \, dx^{1} \wedge dx^{2}.
\ea
\ee
In this case, using (2.5), all $2$-forms (2.6) can be clearly expressed
in the form
\be
\Omega_{s} \equiv \sum_{r=1}^{4} \frac{\partial G_{s \mu}}{\partial \xi^{r}} \,
d \xi^{r} \wedge dx^{\mu} +
\left(\sum_{l=1}^{8} \left(G_{l \nu} \frac{\partial G_{s \mu}}{\partial u^{l}} \right)+
\frac{\partial G_{s \mu}}{\partial x^{\nu}} \right)  dx^{\nu} \wedge dx^{\mu},
\quad s=1, \ldots, 8.
\ee
Let $Y_{\mu}$ be linearly independent vector f\/ields def\/ined on ${\mathbb R}^{14}$
\be
Y_{\mu} = \left(a^{1}_{\mu} \partial_{x^{1}}, a^{2}_{\mu} \partial_{x^{2}},
b^{1}_{\mu} \partial_{\xi^{1}}, \ldots, b^{4}_{\mu} \partial_{\xi^{4}},
c^{1}_{\mu} \partial_{u^{1}}, \ldots, c^{8}_{\mu} \partial_{u^{8}}\right), 
\qquad \mu=1,2
\ee
such that they annihilate systems (2.4) and (2.6),
composed of the $1$-forms $\omega_{s}$ and the
$2$-forms $\Omega_{s}$, respectively, that is
\be
\langle \omega_{s} \lrcorner  Y_{\mu}\rangle =0 , \qquad
\langle \Omega_{s} \lrcorner  Y_{1}, Y_{2} \rangle = 0,  \qquad s=1, \ldots, 8, \quad \mu=1,2
\ee 
at some regular point $(x, \xi, u) \in {\mathbb R}^{14}$. 
The above system is called a system of polar equa\-tions~[13]. 
The set of vector f\/ields
$Y_{\mu}$ satisfying this system depends on a certain 
number~$N$ of free parameters.
In our case, the solution of (2.9) is given by
\[
\ba{l}
\ds Y_{1} = \partial_{x^{1}} + \sum_{r=1}^{4} b_{1}^{r} \partial_{\xi^{r}}+
\xi^{1} \partial_{u_{1}} -p u_{1} \partial_{u_{2}}+p u_{4} \partial_{u_{3}}
+\xi^{4} \partial_{u_{4}}
-[u_{2} \bar{\partial} p - p^{2} u_{1}] \partial_{u_{5}} 
\vspace{3mm}\\
\ds \qquad + [u_{1} \bar{\partial} p + p \xi^{1}] \partial_{u_{6}} -
[ u_{4} \bar{\partial} p + p \xi^{4}] \partial_{u_{7}}+
[u_{3} \bar{\partial} p + p^{2} u_{4}] \partial_{u_{8}},
\ea
\]
and
\be
\ba{l}
\ds Y_{2} = \partial_{x^{2}} + \sum_{r=1}^{4} b_{2}^{r} \partial_{\xi^{r}}
+ p u_{2} \partial_{u_{1}} + \xi^{2} \partial_{u_{2}}
+ \xi^{3} \partial_{u_{3}} -p u_{3} \partial_{u_{4}}
-[u_{2} \partial p + p \xi^{2}] \partial_{u_{5}} 
\vspace{3mm}\\
\ds \qquad + [u_{1} \partial p +
p^{2} u_{2}] \partial_{u_{6}} -[ u_{4} \partial p - p^{2} u_{3}]
\partial_{u_{7}} + [ u_{3} \partial p + p \xi^{3}] \partial_{u_{8}},
\ea
\ee
where
\[
\ba{lll}
b_{1}^{2} = -(u_{1} \partial p + p^{2} u_{2}), & \qquad &
b_{1}^{3} = u_{4} \partial p - p^{2} u_{3},  
\vspace{2mm}\\
b_{2}^{1} = u_{2} \bar{\partial} p - p^{2} u_{1}, & &
b_{2}^{4} = -(u_{3} \bar{\partial} p + p^{2} u_{4}).
\ea
\]
To simplify formulae (2.10) we have used notation (2.3).
Solution (2.10) contains four arbitrary parameters
$b_{1}^{1}$, $b_{1}^{4}$,  $b_{2}^{2}$, $b_{2}^{3}$, hence we have
\be
N=4.
\ee
According to the def\/inition of the f\/irst Cartan character [13], we have
\[
s_{1} = 
\mbox{max\, rank}_{X=(X^{1}, X^{2}) \in {\mathbb R}^{2}}  
\left| \begin{array}{ccc}
\ds \frac{\partial G_{1 \mu}}{\partial \xi^{1}} X^{\mu}, & \cdots &
\ds \frac{\partial G_{1 \mu}}{\partial \xi^{4}} X^{\mu}   \\
\vdots   &       &   \vdots   \\
\ds \frac{\partial G_{8 \mu}}{\partial \xi^{1}} X^{\mu} , &  \cdots &
\ds \frac{\partial G_{8 \mu}}{\partial \xi^{4}} X^{\mu}   
\end{array}   \right|
\]
at a regular point $(x,\xi,u) \in{\mathbb R}^{14}$.
The nonvanishing elements of the $ 8 \times 4$ matrix
$(a_{sr}) =\left(\frac{\partial G_{s \mu}}{\partial \xi^{r}}X^{\mu}\right)$ are
\[
\ba{l}
a_{11} =X^{1}, \qquad a_{22}=X^{2}, \qquad a_{33} =X^{2}, \qquad a_{44}=X^{1},
\qquad a_{51}=u_{2} u_{3}X^{1},
\vspace{2mm}\\
\ds  a_{52}=(p+u_{2}u_{4})X^{2}, \qquad \!
a_{53} = u_{1}u_{2} X^{2}, \qquad \!  a_{54}=  u_{2}^{2} X^{1},
\qquad \! a_{61} =-(p+u_{1}u_{3})X^{1}, \!
\vspace{2mm}\\
a_{62} =-u_{1}u_{4} X^{2}, \qquad  
a_{63} =-u_{2}^{2} X^{2}, \qquad   a_{64} =-u_{1}u_{2} X^{1},
\qquad    a_{71}=u_{3}u_{4}X^{1}, 
\vspace{2mm}\\
\ds a_{72}=u_{4}^{2} X^{2}, \qquad
a_{73}= u_{1}u_{4} X^{2},  \qquad a_{74} =(p+u_{2}u_{4})X^{1},
\qquad  a_{81}=-u_{3}^{2} X^{1}, 
\vspace{2mm}\\
\ds a_{82} =-u_{2}u_{3} X^{1}, \qquad
a_{83} =-(p+u_{1}u_{3}) X^{2}, \qquad a_{84} =-u_{3}u_{4} X^{2},
\ea
\]
since the function $G_{s \mu}$ depends linearly on $\xi$.
Hence, the maximal rank of the matrix $(a_{sr})$~is 
\[
s_{1}=4.
\]
In that case, the second Cartan character is given by
\[
s_{2} = n - s_{1} =0,
\]
where $n=4$ is the number of coordinates $\xi$. Taking into
account the def\/inition~[13] of the Cartan number $Q$, we have
\be
Q= s_{1} + 2 s_{2} = 4.
\ee
Thus, from (2.11) and (2.12), we get
\[
Q=N=4
\]
and, according to Cartan's Theorem, system (2.4) 
is in involution at the regular point $(x_{0}, \xi_{0}, u_{0})$.
So, its general analytic solution exists in some neighbourhood of
this regular point and depends on four arbitrary real analytic
functions of one real variable. \hfill \rule{3mm}{3mm} 

Let us note that, since system (2.4) is equivalent to WE
system (1.2), then Proposition~1 implies the existence of the
general analytic solution of (1.2).
This solution depends on two arbitrary complex analytic functions
of one complex variable and their complex conjugate functions
(since we interpret $z$ and $\bar{z}$ as coordinates on $\mathbb C$
and $\psi_{i}$ and $\bar{\psi}_{i}$ as complex conjugate functions
on $\mathbb C$).

\setcounter{equation}{0}

\section{On the boundary value problem for \\ 
the Weierstrass-Enneper system}

We start by considering the connection between the
structure of certain conserved quantities associated with 
WE system (1.2) and the possibility of the construction of 
some classes of potential solutions.

From the conservation law associated with system (1.2) 
\be
\partial (\psi_{1})^{2} + \bar{\partial} (\psi_{2})^{2} =0,   \qquad
\bar{\partial} (\bar{\psi}_{1})^{2} + \partial (\bar{\psi}_{2})^{2} =0,
\ee
it follows that there exists a potential 
function $g(z, \bar{z}) :{\mathbb C} \rightarrow {\mathbb C}$, such that
the functions $\psi_{i}$ can be expressed in terms of the f\/irst
derivatives of the function $g$
\be
\ba{l}
\ds \psi_{1} = e^{in \pi} (\bar{\partial} g)^{1/2},  \qquad
\psi_{2} = i e^{ik \pi} (\partial g)^{1/2},
\vspace{2mm}\\
\ds \bar{\psi}_{1} = e^{-in \pi} (\partial \bar{g})^{1/2}, \qquad
\bar{\psi}_{2} =-i e^{-ik \pi} ( \bar{\partial} \bar{g})^{1/2},
\qquad n,k \in {\mathbb Z}.
\ea
\ee
Substituting (3.2) into WE system (1.2) one obtains
\be
\ba{l}
\partial \bar{\partial} g = 2i e^{i(k-n)\pi} \left[(\bar{\partial} g)
(\partial \bar{g})^{1/2} (\partial g)^{1/2} + (\bar{\partial} g)^{1/2}
(\partial g)( \bar{\partial} \bar{g})^{1/2}\right],    
\vspace{3mm}\\
\ds \bar{\partial} \partial \bar{g} = -2i e^{-i(k-n) \pi}\left[(\partial \bar{g})
(\bar{\partial} g)^{1/2} (\bar{\partial} \bar{g})^{1/2} +
(\partial \bar{g})^{1/2} (\bar{\partial} \bar{g} )(\partial g)^{1/2}\right].
\ea
\ee
This result can be summarized as follows.

\medskip

\noindent
{\bf Proposition 2.} {\it If a complex-valued
function $g$ of the class $C^{2}$ is a
solution of system~(3.3), then the complex
valued functions $\psi_{i}$ defined by (3.2) are
solutions of WE system~(1.2).}

\medskip

In the next section we show some examples of these types of solutions.

Let us now establish the existence and uniqueness of the
potential solution to the boundary value problem (BVP) for
WE system (1.2). The BVP for this system consists in f\/inding
a class of solutions $\psi_{i}$ in some open bounded simply
connected region $\Omega$ in $\mathbb C$ for prescribed
values of the functions $\psi_{i}$ along the boundary $\partial \Omega$
\be
\ba{l}
\ds \partial \psi_{1} = (|\psi_{1}|^{2} + |\psi_{2}|^{2}) \psi_{2}, \qquad
\bar{\partial} \psi_{2} = - (|\psi_{1}|^{2} + |\psi_{2}|^{2}) \psi_{1},
\qquad \mbox{in} \quad \Omega,
\vspace{3mm}\\
\ds \psi_{1} (z, \bar{z}) = e^{in \pi} (\bar{\partial} g)^{1/2}|_{\partial \Omega},
\qquad \psi_{2} = i e^{i k \pi} (\partial g)^{1/2} |_{\partial \Omega},
  \qquad \mbox{on} \quad \partial \Omega.
\ea
\ee
We show now how a certain class of dif\/ferentiable
solutions of this problem can be obtained with the help
of the conservation law (1.7).
Substituting (3.2) into (1.6), we get 
\be
\partial^{2} g - \frac{\partial g}{\partial \bar{g}} \partial^{2} \bar{g}
+2i e^{i(n-k) \pi} j(z) \left(\frac{\partial g}{\partial \bar{g}}\right)^{1/2} =0.
\ee
The condition for the existence of the entire function
$j(z)$ requires that
\be
\bar{\partial}\left[ \frac{(\partial \bar{g})^{1/2}}{(\partial g)^{1/2}}
\partial^{2} g - \frac{(\partial g)^{1/2}}{(\partial \bar{g})^{1/2}}
\partial^{2} \bar{g}\right] = 0.
\ee
It should be noted that condition (3.6) is identically satisf\/ied
whenever equations~(3.3) hold. This fact simplif\/ies considerably
the process of solving the BVP for the WE system,
since the potential solutions
of WE system (1.2) have to satisfy only conditions (3.3).
Under the above considerations, we can formulate the following.

\medskip

\noindent
{\bf Proposition 3.}  {\it The solution of the 
boundary value problem (3.4) for WE system 
on a simply connected region $\Omega$ exists and is unique,
provided that there exists a $C^{2}$ complex valued function $g$
such that the first order derivatives of 
the function $g$ satisfy equations~(3.3) on $\partial \Omega$.}

\medskip

\noindent
{\bf Proof.} Indeed, if the values of the derivatives 
$\partial g$ and $\bar{\partial} g$ are given on the boundary~$\partial \Omega$, 
such that equations (3.3) hold, then the functions
$\psi_{i}$, def\/ined by (3.2), satisfy WE system~(1.2). This 
fact follows from Proposition~2.
The conservation law (3.5) implies that the current $j$
is an entire function determined by
\[
j(z) = \frac{i}{2} e^{i(k-n) \pi} \frac{(\partial \bar{g})^{1/2}}
{(\partial g)^{1/2}} \left[ \partial^{2} g - \frac{\partial g}{\partial \bar{g}}
\partial^{2} \bar{g} \right] \qquad \mbox{on} \quad \partial \Omega.
\]
Liouville's Theorem ensures that the values of the entire function $j(z)$ are
uniquely def\/ined on the whole simply 
connected region $\Omega \subset \mathbb C$.
This means that there exists a one-to-one correspondence between functions
$\psi_{i}$ prescribed on $\partial \Omega$ and the entire function
$j(z)$. Hence, the values of the solutions $\psi_{i}$ of WE system~(1.2)
at the point $z \in \Omega$ depend only on the values of the
functions $\psi_{i}$ on the boundary $\partial \Omega$. This being so,
the functions $\psi_{i}$ def\/ined by equations (3.2), with the property
that the f\/irst derivatives of the
function $g$ satisfy~(3.3) on $\partial \Omega$, 
are the unique solutions of the BVP (3.4) for WE system 
in the region~$\Omega$. \hfill \rule{3mm}{3mm} 

\setcounter{equation}{0}

\section{Potential solutions of the Weierstrass-Enneper system}

Proposition 2 provides us with a tool for constructing
particular classes of potential solutions to WE system (1.2).
We now present a couple of examples of such solutions.

{\bf 1.} A class of rational solutions of system (3.3) is given by
\be
g(z, \bar{z}) = - \frac{z^{m}}{1 + |z|^{2m}}, \qquad  m \in \mathbb Z.
\ee
Substituting the function $g$ into relations (3.2), one obtains
an explicit solution of WE system (1.2) 
\be
\psi_{1} = e^{i n \pi} m^{1/2} \frac{|z|^{m} }
{1+ |z|^{2m} } z^{(m-1)/2},  \qquad
\psi_{2} = e^{i k \pi} m^{1/2} \frac{z^{(m-1)/2}}{1 + |z|^{2m} }.
\ee
It is interesting to note that this very simple and direct method
yields the same result that was obtained by
a much more complex approach via
the $SU(2)$ sigma model~[14].
For every f\/ixed $m$, solution~(4.2)
belongs to a given topological sector. 
The solutions are double valued for all even $m$.
Each solution (4.2) corresponds to a particular constant mean curvature
surface which is covered $n$ times as $z$ runs over the 
complex plane $\mathbb C$.
This surface is obtained by the parametrization (1.3)
\setcounter{equation}{2}
\renewcommand{\theequation}{\arabic{section}.\arabic{equation}.1}
\be
X_{1} + i X_{2} = 2i z^{-m} \left(\frac{1- |z|^{2m} }{1 + |z|^{2m} }\right),
\ee
\setcounter{equation}{2}
\renewcommand{\theequation}{\arabic{section}.\arabic{equation}.2}
\be
X_{1} - i X_{2} =-2i \bar{z}^{-m} \frac{1 - |z|^{2m} }{1 + |z|^{2m} },
\ee
\setcounter{equation}{2}
\renewcommand{\theequation}{\arabic{section}.\arabic{equation}.3}
\be
X_{3} = \frac{4}{1 + |z|^{2m} }.
\ee
Solving expression (4.3.3) for $|z|^{m}$ in terms of $X_{3}$ and substituting
the result into the expressions (4.3.1) and (4.3.2), one arrives at the
equation of a surface which is obtained by revolving the curve
\[
X_{2}  = 2 \left(2 X_{3} - 1 \right) \left(\frac{X_{3}}{1-X_{3}}\right)^{1/2}
\]
around the $X_{3}$ axis. This surface has a conic point at $(0,0,2)$,
and the corresponding Gaussian curvature is $K=1$.

{\bf 2.} Let us discuss now the 
construction of an algebraic
multi-soliton solution to WE system (1.2). First, we look for a
particular class of rational solutions $g$ of (3.3) 
admitting two simple poles only.
This leads us to the following solution 
\setcounter{equation}{3}
\renewcommand{\theequation}{\arabic{section}.\arabic{equation}}
\be
\ba{l}
\ds 
g(z, \bar{z}) = \frac{(a-b)^{2} (z-a)^{2}}
{(-(2z -a-b)\bar{z} + (z-a)a + (z-b)b)(2z-a-b)} 
\vspace{3mm}\\
\ds \qquad \qquad \qquad + \frac{(a-b)^{2}} {2(2z-a-b)}, \qquad  a,b \in {\mathbb R}.
\ea
\ee 
The corresponding solution of WE system (1.2) takes the form 
\be
\psi_{1} = e^{in \pi} (a-b) \frac{z-a}{|z-a|^{2} + |z-b|^{2}}, \qquad
\psi_{2} = e^{ik \pi} (a-b) \frac{\bar{z} -b}{|z-a|^{2} + |z-b|^{2}}.
\ee
This type of solution is known in the literature [10], and
represents a one-soliton solution.
The associated surface can be computed
from equations (1.3), namely 
\be
\ba{l}
\ds X_{1} + iX_{2}= 2i (a-b)^{2} \left( - \frac{ (\bar{z} -a)^{2}}{\bar{D}}
+\frac{ (z-b)^{2}}{D} \right),
\vspace{3mm}\\
\ds X_{1} - i X_{2} = 2i (a-b)^{2} \left( - \frac{ (\bar{z}-b)^{2}}{\bar{D}}
+ \frac{ (z-a)^{2}}{D}\right),
\vspace{3mm}\\
\ds X_{3} = -2 (a-b)^{2} \left(- \frac{(\bar{z}-a)(\bar{z}-b)}{\bar{D}} 
- \frac{(z-a)(z-b)}{D} \right),
\ea
\ee
where the denominator $D$ is given by
\[
D= ((2z -a-b) \bar{z} -a(z-a) -b(z-b))(2z-a-b), 
\]
and its respective complex conjugate is
\[
\bar{D} = ((2 \bar{z} -a-b) z -a(\bar{z} -a) -b(\bar{z}-b)) (2 \bar{z} -a-b).
\]
Finally, we solve (4.6) for $X_{1}$ and $X_{2}$ and write the
$X_{j}$, $(j=1,2,3)$ in terms of 
$z=x+iy$ and $\bar{z}=x-iy$. This gives us the parametric
forms for the $X_{j}$ in terms of $x$ and $y$. Eliminating
the variables $x$ and $y$ from these equations, one can write
an expression for the surface just in terms of the coordinates $X_{j}$
as follows,
\be
X_{2}^{3} + X_{1}^{2} X_{2} + X_{2} X_{3}^{2} +
4(b-a) \left(X_{1}^{2} + X_{2}^{2} +X_{3}^{2}\right) 
+4 (a -b)^{2} X_{2} = 0.
\ee
The curvature for this surface is given by $K = (a-b)^{-2}$.
Thus, formula (4.7) represents an Enneper type surface.

We now consider a more general case when the solution $g$ of (3.3)
admits arbitrary number of simple poles. Under this assumption we have
\be
\ba{l}
\ds \partial g= -
 \frac{1}{\left(1+  \prod\limits_{j=1}^{N} |\frac{z - a_{j}}{z- b_{j}}|^{2}\right)^{2}}
\left(\sum_{s=1}^{N} \frac{1}{(z-b_{s})} \left(\prod_{j=1 \atop j \neq s}^{N}
\! \frac{(z-a_{j})}{(z -b_{j})} - \prod\limits_{j=1}^{N} \! \frac{(z- a_{j})}{(z-b_{j})}\right)\! \right)\!,
\vspace{3mm}\\
\ds \partial \bar{g} =  
\frac{\prod\limits_{j=1}^{N} \frac{\bar{z} - a_{j}}{\bar{z}-b_{j}}
}{\left(1+  \prod\limits_{j=1}^{N} |\frac{z- a_{j}}{z-b_{j}}|^{2}\right)^{2}} 
\left( \sum_{s=1}^{N} \frac{1}{(z-b_{s})} \left(\prod_{j=1 \atop j \neq s}^{N}
\! \frac{(z-a_{j})}{(z-b_{j})}  - \prod\limits_{j=1}^{N}\! \frac{(z-a_{j})}{(z-b_{j})}\right)\!\right) \! ,
\ a_{j}, b_{j} \in {\mathbb R}.
\ea\hspace{-24.4pt}
\ee
Substituting (4.8) into (3.2) we determine explicitly the 
corresponding form of an algebraic multi-soliton solution of WE system
(1.2)
\be
\ba{l}
\ds \psi_{1} = e^{in \pi} \frac{ \prod\limits_{j=1}^{N} \frac{z-a_{j}}{z- b_{j}}}
{1 +  \prod\limits_{j=1}^{N} |\frac{z-a_{j}}{z-b_{j}}|^{2}}
\left( \sum_{s=1}^{N}  \frac{1}{(\bar{z} - b_{s})} \left(
 \prod_{j=1 \atop j \neq s}^{N} \frac{(\bar{z} -a_{j})}{(\bar{z} -b_{j})}
- \prod_{j=1}^{N} \frac{(\bar{z} -a_{j})}{(\bar{z}-b_{j})}\right)\right)^{1/2},  
\vspace{3mm}\\
\ds \psi_{2} = 
 \frac{e^{ik \pi}}{1+  \prod\limits _{j=1}^{N} |\frac{z- a_{j}}{z-b_{j}}|^{2}}
\left( \sum_{s=1}^{N} \frac{1}{(z-b_{s})} \left( \prod_{j=1 \atop j \neq s}^{N}
\frac{(z-a_{j})}{(z-b_{j})} - \prod_{j=1}^{N} \frac{(z-a_{j})}{(z-b_{j})}\right)\right)^{1/2}.
\ea
\ee 
Note that the topological charge (1.5) for each of the instanton 
solutions (4.5) entering into the superposition corresponds to an integer $I= e^{i n \pi} N$.
It is interesting to note that the constant mean curvature surface corresponding to (4.9)
is also determined by (4.7).

\setcounter{equation}{0}

\section{Harmonic solutions of the Weierstrass-Enneper system}

We discuss now the existence of a class of harmonic solutions 
to the WE system (1.2) which can be obtained by applying certain composition
transformations. 

\medskip

\noindent
{\bf Proposition 4.} {\it Suppose that the complex valued
functions $f_{i}$ and $\bar{f}_{i}$ are solutions of
the following system of differential equations,
\be
\ba{l}
\ds p f_{1}''(v) f_{1} (v) f_{2} (v) 
\vspace{3mm}\\
\ds \qquad +\left[- |f_{1}|^{2} (f_{1}' (v))^{2} +
|f_{2}|^{2} \bar{f}'_{2} (\bar{v}) \frac{f_{1}'(v) f_{2}'(v)}{\bar{f}_{1}'  
(\bar{v})}\right] f_{2}(v) - p f_{1}(v) f_{1}'(v) f_{2}'(v) = 0,  
\vspace{3mm}\\
\ds p f_{2}'' (v) f_{1} (v) f_{2} (v) 
\vspace{3mm}\\
\ds \qquad + \left[|f_{1}|^{2} \frac{f_{1}'(v) f_{2}'(v)}
{\bar{f}_{2}' (\bar{v})} \bar{f}_{1}'(\bar{v}) - |f_{2}|^{2}
(f_{2}'(v))^{2}\right] f_{1} (v) - p f_{2} (v) f_{1}'(v) f_{2}' (v) = 0,
\ea
\ee
with respect to the relevant variables $v$ and $\bar{v}$.
Then the compositions of the functions $f_{i}$ with any
harmonic function $v= h(z,\bar{z})$, defined in a simply connected
region $\Omega$, 
\[
\psi_{i} = f_{i} (h(z, \bar{z})),  \qquad  i=1,2
\]
and their respective complex conjugates 
\[
\bar{\psi}_{i} = \bar{f}_{i} ( \bar{h}(\bar{z}, z)),
\qquad i=1,2
\]
constitute solutions of the WE system (1.2).}

\medskip

\noindent
{\bf Proof.} Substituting the composed functions 
$\psi_{i}$ and $\bar{\psi}_{i}$ into (1.2), one obtains 
\be
\ba{l}
f_{1}' (v) \, \partial h (z, \bar{z}) = p f_{2} ( h(z, \bar{z})), \qquad
f_{2}' (v) \, \bar{\partial} h(z, \bar{z} ) = - p f_{1} ( h(z, \bar{z})), 
\vspace{2mm}\\
p =|f_{1}|^{2} + |f_{2}|^{2}. 
\ea
\ee
Dif\/ferentiating the f\/irst equation
in (5.2) with respect to $\bar{\partial}$ and making use of
the equation $\bar{\partial} \partial h = 0$, one obtains 
\be
f_{1}'' (v) \bar{\partial} h \partial h = (\bar{\partial} p) f_{2} +
p f_{2}' (v) \, \bar{\partial} h.
\ee
Similarly, dif\/ferentiating the second equation of (5.2) with respect to
$\partial$, and taking into account the relation
$\bar{\partial} \partial h = 0$, one gets
\be
f_{2}'' (v) \, \partial h \, \bar{\partial} h = - (\partial p) f_{1}
- p f_{1}' (v) \, \partial h.
\ee
Dif\/ferentiating $p$ and using (1.2), we have
\be
\partial p = \psi_{1} \partial \bar{\psi}_{1} + \bar{\psi}_{2}
\partial \psi_{2}, \qquad 
\bar{\partial} p = \bar{\psi}_{1} \bar{\partial} \psi_{1} +
\psi_{2} \bar{\partial} \bar{\psi}_{2}.
\ee
Solving (5.2) for the derivatives $\partial h$, $\bar{\partial} h$ 
and substituting these derivatives
and the expressions (5.5) into equations
(5.3) and (5.4), one obtains
the system of dif\/ferential equations~(5.1). Thus, 
equations~(5.1) are equivalent to (1.2) whenever
$v$ is a harmonic function.~\hfill~\rule{3mm}{3mm} 

Let us consider a simple example to illustrate 
Proposition 4.
A special class of exponential solutions of (5.1) has the form
\[
f_{1} = -i e^{i v} , \qquad 
f_{2} = a e^{i v},   \qquad a \in {\mathbb C}.
\]
If we choose a specif\/ic harmonic form of the function
$v = q (z^{2} + \bar{z}^{2})$, $q \in \mathbb R$,
then the compositions of the functions $f_{i}$ and $v$ give
particular solutions of the WE system
\be
\psi_{1} = -i e^{i q (z^{2} -  \bar{z}^{2} )}, \qquad
\psi_{2} = a e^{i q ( z^{2} -  \bar{z}^{2} )}, \qquad |a|^{2} = 1.
\ee 
The corresponding constant mean curvature surface is determined by
relations (1.3)
\[
\ba{l}
\ds X_{1} +i X_{2} = -i \pi^{1/2} \left(\frac{\mbox{erf}\,(\xi z)}{\exp(2iq \bar{z}^{2}) \xi}
+ \frac{\bar{a}^{2} \exp( 2 i q z^{2}) \, \mbox{erf}\, (\eta \bar{z})}{\eta}\right),
\vspace{3mm}\\
\ds X_1 - i X_{2} = i \pi^{1/2}\left( \frac{a^{2} \mbox{erf}\, ( \xi z)}{\exp(2i q \bar{z}^{2})
\xi} + \frac{\exp(2 i q z^{2}) \, \mbox{erf}\, ( \eta \bar{z})}{\eta}\right),
\vspace{3mm}\\
\ds X_{3} =  - i \pi^{1/2}\left( \frac{ \mbox{erf}\, (\xi z)}{\exp(2 i q \bar{z}^{2}) \xi}
+ \frac{\bar{a} \exp(2 i q z^{2}) \, \mbox{erf}\, ( \eta \bar{z})}{\eta}\right),
\ea
\]
where erf is the error function, $\xi =(-2 i q)^{1/2}$ and
$\eta = (2 i q)^{1/2}$.
The elimination of the quantities 
$\mbox{erf}\,( \xi z)$ and $\exp ( 2i q z^{2})$ 
from the above expressions leads to the formula which
represents a constant mean curvature surface 
describing a catenoide
\be
4\left(1 -  a_{r}^{2}\right) X_{1}^{2} + 4 a_{r}^{2} X_{2}^{2}
+ 8 a_{r} \left(1 - a_{r}^{2}\right)^{1/2} X_{1} X_{2} + 4 X_{3}^{2} = 0, \qquad
a_{r} = Re \, a.
\ee
The Gaussian curvature is $K=1-a_{r}^{2}$.

\setcounter{equation}{0}

\section{Reduction of the Weierstrass-Enneper system \\ to a linear system}

Now we discuss the case when the WE system is subjected to a
single dif\/ferential constraint. This allows us to reduce this system to
a system of linear coupled PDEs.

\medskip

\noindent
{\bf Proposition 5.} {\it The overdetermined system composed of the 
WE system (1.2) and the first order differential constraint
\be
\psi_{1} \partial \bar{\psi}_{1} - \epsilon \bar{\psi}_{1}
\bar{\partial} \psi_{1} + \bar{\psi}_{2} \partial \psi_{2} -
\epsilon \psi_{2} \bar{\partial} \bar{\psi}_{2} = 0,  \qquad
\epsilon = \pm 1
\ee
is consistent if the conditions 
\be
| \psi_{1}|^{2} + |\psi_{2}|^{2} = p(z + \epsilon \bar{z}),
\ee
and 
\be
\ddot{p} - \frac{\dot{p}^{2}}{p} + \epsilon \frac{A}{p} 
- \epsilon p^{3} = 0,    \qquad A \in {\mathbb R}^{+},
\ee
hold. The general analytic solution 
of the above overdetermined system
depends on one arbitrary complex analytic function
of one complex variable and on its complex conjugate function.}

\medskip

\noindent
{\bf Proof.} Indeed, from equations (5.5), taking into account the dif\/ferential
constraint (6.1), one obtains
\be
(\partial - \epsilon \bar{\partial} ) p = \psi_{1} \partial \bar{\psi}_{1} +
\bar{\psi}_{2} \partial \psi_{2}  - \epsilon \bar{\psi}_{1}
\bar{\partial} \psi_{1} - \epsilon \psi_{2} \bar{\partial}
\bar{\psi}_{2} = 0.
\ee
This means that, under the assumption (6.1), the quantity $p$ is 
a real valued function of the argument $s = (z + \epsilon \bar{z})$. 
Hence, $p$ is a conserved quantity. Therefore, condition (6.2) holds.
Dif\/ferentiating equation (6.2) with respect to $z$ and $\bar{z}$, one obtains
\be
\psi_{1} \partial \bar{\psi}_{1} + \bar{\psi}_{2} \partial \psi_{2}
= \dot{p},  \qquad
\bar{\psi}_{1} \bar{\partial} \psi_{1} + \psi_{2}
\bar{\partial} \bar{\psi}_{2} = \epsilon \dot{p},
\ee
where we introduced the notation $\dot{p} = dp/d s$. 
Solving equation (6.5) with respect to $\partial \psi_{2}$ and next
substituting this term into (1.6), one obtains
\be
\partial \bar{\psi}_{1}  = \frac{\bar{\psi}_{1}}{p}
\left(\dot{p} - \frac{\bar{\psi}_{2}}{\bar{\psi}_{1}} j(z)\right).
\ee
The complex conjugate of equation (6.6) is given by
\[
\bar{\partial} \psi_{1} = \frac{\psi_{1}}{p}
\left(\epsilon \dot{p} - \frac{\psi_{2}}{\psi_{1}} \bar{j} (\bar{z})\right).
\]
A similar analysis can be performed for $\partial \bar{\psi}_{1}$,
in order to determine derivatives of $\psi_{2}$ in terms of
$\psi_{i}$, $j$ and $p$. As a result, 
one obtains from WE system (1.2) the following
system of equations
\be
\ba{lll}
\ds \partial \psi_{1} =p \psi_{2},  & \qquad &
 \ds \partial \psi_{2} = \frac{\psi_{2}}{p} \left( \dot{p} + \frac{\psi_{1}}{\psi_{2}} j(z)\right),   
\vspace{3mm}\\
\ds  \bar{\partial} \psi_{1} = \frac{\psi_{1}}{p} \left(\epsilon \dot{p}
- \frac{\psi_{2}}{\psi_{1}} \bar{j}(\bar{z})\right),   &  &
\ds  \bar{\partial} \psi_{2} = - p \psi_{1}, 
\vspace{3mm}\\
\ds  \partial \bar{\psi}_{1} = \frac{\bar{\psi}_{1}}{p}
\left(\dot{p} - \frac{\bar{\psi}_{2}}{\bar{\psi}_{1}} j(z)\right),  & &
\ds  \partial \bar{\psi}_{2} = - p \bar{\psi}_{1},  
\vspace{3mm}\\
\ds  \bar{\partial} \bar{\psi}_{1} = p \bar{\psi}_{2},  & &
\ds  \bar{\partial} \bar{\psi}_{2} =
\frac{\bar{\psi}_{2}}{p} \left(\epsilon \dot{p} + \frac{\bar{\psi}_{1}}
{\bar{\psi}_{2}} \bar{j}(\bar{z})\right).
\ea
\ee
The compatibility conditions for (6.7) 
require $|j|^{2} = A \in {\mathbb R}^{+}$ and the current $j$ is constant.
This gives a dif\/ferential equation for $p$ of the form (6.3). 
Expressing~(6.7) in the language of dif\/ferential one-forms
and making use of Propositions~1, one can easily show that,
if the compatibility conditions (6.2) and (6.3) are satisf\/ied, then
the general analytic solution of~(6.7) with $j$ constant exists and depends
on one arbitrary analytic complex function
of one complex variable.\hfill \rule{3mm}{3mm} 

Now, let us discuss some classes of solutions of the dif\/ferential
equation (6.3) which allow one to reduce the WE system to the
coupled linear PDEs. Note that equation (6.3) is of Painlev\'{e} type, PXII, having
only poles for moveable singularities.
The f\/irst integral is given by
\be
\dot{p}(s)^{2} = (\epsilon p^{4} + K p^{2} + \epsilon A),
\ee
where $K$ is an arbitrary real constant. 
The forms of the real solutions for $p$ depend on the relationships
between the roots of the right-hand side of the ODE (6.8).
They lead to the following cases.

(i) Elementary solutions, such as constant, algebraic with one or two
simple poles, trigonometric and hyperbolic solutions.

(ii) Doubly periodic solutions which can be expressed in terms of
the Jacobi elliptic functions $\mbox{sn}$ and $\mbox{cn}$. The moduli $k$ of
the elliptic functions are chosen in such a way that $0 < k^{2} <1$.
This fact ensures that the elliptic solutions possess one real and 
one purely imaginary period. Consequently, for real argument $s$, we have
$-1 \leq \mbox{sn}\, (s,k) \leq 1$, and
$-1 \leq \mbox{cn}\, (s,k) \leq 1$.

The dif\/ferent classes of solutions of (6.8)
are summarized in Tables 1 and 2. They lead us to twelve dif\/ferent types
of solutions of the linear coupled WE system for which the compatibility
conditions are satisf\/ied identically.

Finally, let us discuss the case when we add multiple dif\/ferential
constraints compatible with our basic WE system (1.2).
We show that if the conditions (6.1) for $\epsilon =+1$ and $\epsilon =-1$ are
simultaneously satisf\/ied
\be
\psi_{1} \partial \bar{\psi}_{1} + \bar{\psi}_{2} \partial \psi_{2} =0,
\qquad
\bar{\psi}_{1} \bar{\partial} \psi_{1} + \psi_{2} \bar{\partial}
\bar{\psi}_{2} =0
\ee
then WE system (1.2) can be reduced to a linear decoupled system.

\newpage

\begin{landscape}
\begin{table}
\small 

\caption{Finite real elementary and elliptic solutions $p=p(s-s_{0})$ with
$\dot{p}^{2} = \epsilon p^{4} +K p^{2} + \epsilon A$,  where 
$A \geq 0$, $K \in \mathbb R$, $w \equiv \sqrt{K^{2} -4A}$,
$a \equiv (K \pm w)^{1/2}/\sqrt{2}$, $B \equiv 1+ \sqrt{2}$,
$\epsilon, \epsilon_{1} = \pm 1$
and $s-s_{0} =z +\epsilon \bar{z} -s_{0}$.}

\vspace{3mm}

\begin{tabular}{|c|c|c|c|c|c|}  
 \hline
$No$  &  $\epsilon$ & $p(s-s_{0})$  & $k$ & $g$  & $Range$   \\  \hline
\mbox{} & \mbox{}  &  \mbox{}  &  \mbox{}  &  \mbox{}  &  \mbox{}  \\[-3mm]
$1$   &  $+1$       &  $ 
\ds \epsilon_{1} A^{1/4} 
 \left(\frac{1+ \mbox{cn}\left(2 A^{1/4} (s-s_{0}), \frac{1}{2}\right)}
{1-\mbox{cn} \left(2 A^{1/4} (s-s_{0}), \frac{1}{2}\right)}\right)^{1/2} $ & 
$\ds  \frac{1}{\sqrt{2}}$ & --- & $K=0$, \quad $0 \leq p < \infty$ \\[6mm]
$2$ & $+1$ &  $ \ds \epsilon_{1}A^{1/4} \frac{B\, \mbox{cn} \left(\frac{A^{1/4}}{g} (s-s_{0}),k
\right) - \mbox{sn} \left(\frac{A^{1/4}}{g}(s-s_{0}),k\right)}{B\, \mbox{cn}
\left(\frac{A^{1/4}}{g}(s-s_{0}),k\right) + \mbox{sn} \left(\frac{A^{1/4}}{g} (s-s_{0}),k\right)}$ &
$ 2\left( 3\sqrt{2} -4\right)^{1/2}$ & $2-\sqrt{2}$ & $K=0$, \quad $0 \leq p \leq 1$   \\[6mm]
$3$ & $-1$ & $ \ds \frac{\epsilon_{1} (2A)^{1/2}}
{\left[(K+w)(1-k^{2} \mbox{sn}^{2} \left(\frac{(s-s_{0})}{g},k\right)\right]^{1/2}}$ &
$ \ds \left(\frac{2 w}{K+w}\right)^{1/2}$ &
$ \ds \frac{\sqrt{2}}{(K + w)^{1/2}}$ &
$ K>0$,\quad $K^{2}-4A >0, $  \\ 
\mbox{} & \mbox{} & \mbox{} & \mbox{} & \mbox{} & $\ds \frac{1}{\sqrt{2}}
(K-w)^{1/2} \leq p < \frac{1}{\sqrt{2}}(K+w)^{1/2}$  \\[6mm]
$4$ & $-1$ & $\ds \epsilon_{1}  
\left(\frac{1}{2}(K+w)-w \, \mbox{sn} \left(\frac{s-s_{0}}{g},k\right)\right)^{1/2}$ &
$ \ds \left(\frac{2 w}{K+w}\right)^{1/2}$ & 
$ \ds \left(\frac{2}{K+w}\right)^{1/2}$ &
$  K>0$, \quad $K^{2} -4 A >0$  \\
\mbox{} & \mbox{} & \mbox{} & \mbox{} & \mbox{} & $\ds \frac{1}{\sqrt{2}}
(K-w)^{1/2} \leq p < \frac{1}{\sqrt{2}}(K+w)^{1/2}$   \\[6mm]
$5$ & $+ 1$ & $ \ds \epsilon_{1} \frac{(|K| \pm w)^{1/2}}{\sqrt{2} 
\, \mbox{sn} \left(\frac{s-s_{0}}{g},k\right)}$ &
$  \ds \left(\frac{|K| \mp w}{|K| \pm w}\right)^{1/2}$ &
$  \ds \frac{\sqrt{2}}{(|K| \pm w)^{1/2}}$
 & $ K<0$, \quad $K^{2} -4A>0,$  \\
\mbox{} & \mbox{} & \mbox{} & \mbox{} & \mbox{} & $\ds \frac{1}{\sqrt{2}}
(|K| \pm w)^{1/2} \leq p < \infty$  \\[6mm] 
$6$ & $+1$ & $\ds  \frac{K}{\cosh\left(\sqrt{K} (s-s_{0})\right)}$ & --- & --- &
$A=0$, $K>0$, \quad $0<p < K$  \\[6mm]  
  \hline
\end{tabular}
\end{table}
\end{landscape}

\newpage

\begin{landscape}
\small

\begin{table}
\caption{Singular real elementary and elliptic solutions $p=p(s-s_{0})$
with $\dot{p}^{2}=\epsilon p^{4} +K p^{2} + \epsilon A$, where 
$A \geq0$, $K \in \mathbb R$, $\epsilon,\epsilon_{1}=\pm 1$, $w \equiv \sqrt
{K^{2} - 4A}$ and $s-s_{0} = z + \epsilon \bar{z} - s_{0}$.}

\vspace{3mm}

\begin{tabular}{|c|c|c|c|c|c|}
 \hline
$No$  & $\epsilon$ & $p(s-s_{0})$  &  $k$  &  $g$  & $Range$ \\ \hline
\mbox{} & \mbox{} & \mbox{} &  \mbox{}  &  \mbox{}  &  \mbox{}  \\[-3mm]
$1$  & $+1$ & $\ds \epsilon_{1}  \frac{(2 A)^{1/2}}{K \pm w} \, tn
\left(\frac{s-s_{0}}{g},k\right)$ & 
$ \ds \frac{(2(K^{2}-4 A \pm K w))^{1/2}}{K \pm w}$ &
$ \ds \left(\frac{1}{2}(K \pm w)\right)^{-1/2}$ & $K>0$, \quad $K^{2}-4A>0,$  \\
\mbox{} & \mbox{} & \mbox{} & \mbox{} & \mbox{} & $A^{1/2}<(K \pm w)/2$, $0<p$  \\[6mm]
$2$ & $+1$ & $\ds \epsilon_{1}  \frac{\left[|K| \pm w+(|K| \mp w) \, \mbox{sn}^{2}
\left( \frac{s-s_{0}}{g},k\right)\right]^{1/2}}
{\sqrt{2} \, \mbox{cn} \left( \frac{s-s_{0}}{g},k\right)}\! $ & 
$\ds \left( \frac{|K| \mp w}{|K| \pm w}\right)^{1/2}$ &
$ \ds \frac{\sqrt{2}}{(|K| \pm w)^{1/2}}$ &
$K<0$, \quad $K^{2}-4A>0,$  \\
\mbox{} & \mbox{} & \mbox{} & \mbox{} & \mbox{} &
$ \ds \frac{1}{\sqrt{2}}(|K| \pm w) < p$  \\[6mm]
$3$ & $+1$ & $a \cdot \mbox{tn}\, (a(s-s_{0}),k)$ & 
$\ds  \frac{(2(K^{2}-4A \pm Kw))^{1/2}}
{K \pm w}$ & $\ds  \frac{\epsilon_{1} \sqrt{2}}{(K \pm w)^{1/2}}$ &
$K>0$, \quad $\ds A^{1/2} < \frac{1}{2} (K \pm w)$  \\[6mm]
$4$ &  $+1$ & $ \ds \frac{K}{\sinh(\sqrt{K}(s-s_{0}))}$ &
--- & --- & $A=0$, \quad $K>0$   \\[6mm]
$5$ &  $+1$ & $ \ds \sqrt{|K|} \sec  \frac{\sqrt{|K|}}{2}(s-s_{0})$ &
--- & --- & $A=0$, \quad $K<0 $, \ $\sqrt{|K|} \leq p$  \\[6mm]
$6$ & $+1$  & $ \ds \frac{1}{s_{0} - s}$ &  ---  &  ---  &  $A=0$, \quad $K=0$,
\quad $0<p$   \\[6mm]
 \hline
\end{tabular}
\end{table}
\end{landscape}

\newpage

\noindent
{\bf Proposition 6.} {\it If the functions $\psi_{1}$ and $\psi_{2}$
satisfying WE system (1.2) are subjected to differential constraints (6.9), then
the WE system is reduced to a linear decoupled system of second
order
\be
\bar{\partial} \partial \psi_{i} + p_{0}^{2} \psi_{i} =0,  \qquad
\partial \bar{\partial} \bar{\psi}_{i} + p_{0}^{2} \bar{\psi}_{i}=0,
\qquad i=1,2,
\ee
where
\be
|\psi_{1}|^{2} + |\psi_{2}|^{2} = p_{0} \in \mathbb R^{+}.
\ee}

\noindent
{\bf Proof.} In fact, taking into account dif\/ferential constraints
(6.9) and the derivatives
of $p$ given by (6.4), we obtain that $p$ is a real positive constant.
This means that the overdetermined system composed of WE system (1.2)
and dif\/ferential constraints (6.9) admits a conserved quantity (6.11). Hence, the
WE system (1.2) can be decoupled into the second order system (6.10).
\hfill \rule{3mm}{3mm} 

A simple example of the solution of (1.2) constructed with the
use of dif\/ferential constraints (6.9) was presented in~[14].
This solution has the form of the plane wave or so called vacuum solution
\be
\psi_{1} = A e^{i(hz+k \bar{z})}, \qquad
\psi_{2} = i \frac{A}{k} e^{i(hz + k \bar{z})}, \qquad hk = p_{0},
\qquad h, k \in {\mathbb R}, \quad A \in {\mathbb C}.
\ee
Proposition 6 implies that, due to the linearity of equations
(6.10), a more general class of solutions can be constructed. 
Namely, the linear superposition of plane waves (6.12) gives
\be
\ba{l}
\psi_{1} = A_{1} e^{\alpha_{1} (z + \bar{z})} + A_{2} e^{\alpha_{2}
(z - \bar{z})}, \qquad
\psi_{2} = B_{1} e^{\alpha_{1} (z + \bar{z})} + B_{2} e^{\alpha_{2}
(z - \bar{z})}, 
\vspace{2mm}\\
\alpha_{i} \in \mathbb R, \qquad A_{i}, B_{i} \in \mathbb C, \qquad i=1,2,
\ea
\ee
where $p_{0}= |A_{1}|^{2} + |A_{2}|^{2} +|B_{1}|^{2} +|B_{2}|^{2}$ and
$A_{1} \bar{A}_{2} + B_{1} \bar{B}_{2} = 0$, 
$\bar{A}_{1} A_{2} + \bar{B}_{1} B_{2} = 0$, $B_{1}= \alpha_{1} A_{1}/p_{0}$
and $B_{2} = \alpha_{2} A_{2} / p_{0}$,
and where $\alpha_{1} = \pm i p_{0}$ and $\alpha_{2} = \pm p_{0}$.
Next, substituting~(6.13) into system
(1.3), we obtain a set of equations which determine
a constant mean curvature surface in the parametric form
\[
\ba{l}
\ds X_{1} + i X_{2} = 2i \left(- \frac{\bar{A}_{1}^{2}}{\alpha_{1}} u^{2}
- \frac{\bar{A}_{2}^{2}}{\alpha_{2}} v^{2} -2 
\left(\frac{1}{\alpha_{1} + \alpha_{2}} + \frac{i}{\alpha_{1} - \alpha_{2}}\right)
\bar{A}_{1} \bar{A}_{2} uv \right),
\vspace{3mm}\\
\ds X_{1} - i X_{2} = 2i \left( - \frac{A_{1}^{2}}{\alpha_{1} u^{2}} +
\frac{A_{2}^{2}}{\alpha_{2} v^{2}} + 2 \left(\frac{i}{\alpha_{1}+ \alpha_{2}}
- \frac{1}{\alpha_{1} - \alpha_{2}}\right) A_{1} A_{2} (uv)^{-1}\right),
\vspace{3mm}\\
\ds X_{3} = \frac{1}{2} \left(\left(i |A_{1}|^{2} + |A_{2}|^{2}\right) 
\left(\frac{\ln (u)}{\alpha_{1}} -\frac{\ln ( v)}{\alpha_{2}}\right)
 - \frac{\bar{A}_{1} A_{2} u}
{(\alpha_{1} - \alpha_{2}) v} + i \frac{A_{1} \bar{A}_{2} v}
{(\alpha_{1}-\alpha_{2}) u}\right)
\vspace{3mm}\\
\ds \qquad + \frac{1}{2} \left(\left(-i|A_{1}|^{2} + |A_{2}|^{2}\right) \left(\frac{\ln(u)}{\alpha_{1}} +
\frac{\ln(v)}{\alpha_{2}}\right)
+ i \frac{\bar{A}_{1} A_{2} u}{(\alpha_{1} + \alpha_{2})v} +
\frac{A_{1} \bar{A}_{2} v}{(\alpha_{1} + \alpha_{2}) u} \right),
\ea
\]
in terms of
$u= \exp(-\alpha_{1}(z+ \bar{z}))$, and $v= \exp(-\alpha_{2}(z - \bar{z}))$.
The Gaussian curvature is $K=1$.

\setcounter{equation}{0}

\section{Separation of variables}

Now let us discuss the separation of variables admitted
by WE system (1.2) which enables us to construct the family of solitonlike
solutions. The methodological 
approach assumed in this section is based on the generalized method of
separation of variables developed in~[17]. We are looking for a special
class of solutions of WE system (1.2) of the form
\be
\psi_{i} (z, \bar{z}) = \varphi_{i} ( X \cdot Y),    \qquad i=1,2,
\ee
where $X=X(z)$, $Y=Y(\bar{z})$. Its respective complex
conjugate is given by
\be
\bar{\psi}_{i} (z, \bar{z}) = \bar{\varphi}_{i} 
(\bar{X} \cdot \bar{Y}),     
\ee
where $\bar{X} = \bar{X}(\bar{z})$, $\bar{Y}=\bar{Y}(z)$.
We assume the existence of two complex scalar functions $\xi$
and $\eta$ of $z$ and $\bar{z}$, respectively, such that the dif\/ferential
equations
\be
\frac{d X}{dz} = \xi(X), \qquad
\frac{d Y}{d \bar{z}} = \eta(Y),
\qquad 
\frac{d \bar{X}}{d \bar{z}} = \bar{\xi} (\bar{X}), \qquad
\frac{d \bar{Y}}{dz} = \bar{\eta} (\bar{Y}),
\ee
hold. This means that the complex functions $X$ and $Y$ are locally
piecewise monotonic functions. Substituting (7.2) into WE system (1.2) and 
taking into account (7.3), we obtain a system of dif\/ferential equations
\renewcommand{\theequation}{\arabic{section}.\arabic{equation}$i$}
\setcounter{equation}{3}
\be
\dot{\varphi}_{1} \xi Y = \left( |\varphi_{1}|^{2} + |\varphi_{2}|^{2}\right)
\varphi_{2},  \qquad
\dot{\bar{\varphi}}_{1} \bar{\xi} \bar{Y} =
\left(|\varphi_{1}|^{2} + |\varphi_{2}|^{2}\right) \bar{\varphi}_{2},
\ee
\renewcommand{\theequation}{\arabic{section}.\arabic{equation}$ii$}
\setcounter{equation}{3}
\be
\dot{\varphi}_{2} \eta X = -\left(|\varphi_{1}|^{2} + |\varphi_{2}|^{2}\right)
\varphi_{1},   \qquad
\dot{\bar{\varphi}}_{2} \bar{\eta} \bar{X} =- 
\left(|\varphi_{1}|^{2} + |\varphi_{2}|^{2}\right) \bar{\varphi}_{1}.
\ee
Let us introduce two dif\/ferential operators
\renewcommand{\theequation}{\arabic{section}.\arabic{equation}}
\setcounter{equation}{4}
\be
A= X \partial_{X} - Y \partial_{Y},   \qquad
\bar{A} = \bar{X} \partial_{\bar{X}} - \bar{Y} \partial_{\bar{Y}},
\ee
which are annihilators of any complex function of $s = X \cdot Y$ and
$\bar{s} = \bar{X} \cdot \bar{Y}$, respectively. These operators
commute
\be
[ A, \bar{A} ] =0.
\ee
We operate with the operators $A$ and $\bar{A}$ on equations
$(7.4i)$ and $(7.4ii)$, respectively. Taking into account (7.3), one obtains
\be
\ba{l}
\ds \dot{\varphi}_{1} (XY \xi' - Y \xi) =0, \qquad
\dot{\bar{\varphi}}_{1} ( \bar{X} \bar{Y} \xi' - \bar{Y} \bar{\xi}) =0,
\vspace{2mm}\\
\dot{\varphi}_{2} ( X \eta - XY \eta') =0,  \qquad
\dot{\bar{\varphi}}_{2} (\bar{X} \bar{\eta} - \bar{X} \bar{Y} \bar{\eta}')=0,
\ea
\ee
where dots or primes mean the derivatives of the respective functions 
with respect to their own arguments. 

Let us consider separately
two cases, namely the case in which $\dot{\varphi}_{i}$ does not vanish
anywhere and the case in which $\dot{\varphi}$ is identically equally to 
zero. It is easy to show that in the case when $\dot{\varphi}_{i} =0$,
equations (7.4) do not admit separable solutions, since
\be
|\varphi_{1} |^{2} + | \varphi_{2} |^{2} = 0
\ee
holds.

In the second case, when
$\dot{\varphi}_{i} \neq 0$, equations (7.7) can be integrated and
their f\/irst integrals are
\be
X = \gamma e^{\alpha z},   \qquad   \bar{X} = \bar{\gamma} e^{\bar{\alpha}\bar{z}}, 
\qquad 
Y= \kappa e^{\beta \bar{z}}, \qquad  \bar{Y} = \bar{\kappa} e^{\bar{\beta} z},
\ee
where $\alpha$, $\beta$, $\gamma$ and $\kappa$ are arbitrary complex constants.
Substituting (7.9) into (7.4), we obtain a nonlinear system of ODEs
\be
\frac{d \varphi_{1}}{ds} = \frac{p}{\alpha s} \varphi_{2}, \qquad
\frac{d \bar{\varphi}_{1}}{d \bar{s}} = \frac{p}{\bar{\alpha} \bar{s}}
\bar{\varphi}_{2},
\qquad 
\frac{d \varphi_{2}}{ds} = - \frac{p}{\beta s} \varphi_{1},   \qquad
\frac{d \bar{\varphi}_{2}}{d \bar{s}} = - \frac{p}{\bar{\beta} \bar{s}}
\bar{\varphi}_{1}.
\ee
System (7.10) can be integrated using the 
condition (1.7) for the conservation of current. 
Taking into account (1.2) and (7.9), we obtain
\be
\alpha^{2} s \left[ \frac{\ddot{\varphi}_{2}(s)}{\varphi_{2}(s)} s
+ \frac{\dot{\varphi}_{2}(s)}{\varphi_{2}(s)}\right] -
\bar{\beta}^{2} \bar{s} \left[ \frac{\ddot{\bar{\varphi}}_{1}(\bar{s})}
{\bar{\varphi}_{1}(\bar{s})} \bar{s} 
 + \frac{\dot{\bar{\varphi}}_{1}(\bar{s})}{\bar{\varphi}_{1}
(\bar{s})} \right] =0.
\ee
The variables $s$ and $\bar{s}$ are separable if the Cauchy-Euler
dif\/ferential equations
\be
s^{2} \ddot{\varphi}_{2} + s \dot{\varphi}_{2} - \frac{\mu}{\alpha^{2}}
\varphi_{2} =0,  \qquad
\bar{s}^{2} \ddot{\bar{\varphi}}_{1} + \bar{s} \dot{\bar{\varphi}}_{1}
- \frac{\mu}{\bar{\beta}^{2}} \bar{\varphi}_{1} =0
\ee
hold, where $\mu$ is a complex separation constant. After the 
integration of (7.12), we obtain the solutions
\be
\ba{l}
\varphi_{1} = d_{1} s^{\bar{q}} + d_{2} s^{-\bar{q}}, \qquad
\varphi_{2} = c_{1} s^{r} + c_{2} s^{-r},
\vspace{2mm}\\
\bar{\varphi}_{1} = \bar{d}_{1} \bar{s}^{q} + \bar{d}_{2}
\bar{s}^{-q},   \qquad
\bar{\varphi}_{2} = \bar{c}_{1} \bar{s}^{\bar{r}}+ \bar{c}_{2}
\bar{s}^{-\bar{r}}.
\ea
\ee
Here, $q,\bar{q}$, $r$, $\bar{r}$ and $c_{i}$, $d_{i}$, $\bar{c}_{i}$,
$\bar{d}_{i}$, for $i=1,2$ are arbitrary complex constants of
integration. Substituting (7.13) into (7.10), we obtain the equations
\be
\ba{l}
\ds -\alpha \bar{q} (d_{1} s^{\bar{q}} - d_{2} s^{-\bar{q}})
+ (c_{1} s^{r} +c_{2} s^{-r}) [ |d_{1}|^{2} s^{\bar{q}} \bar{s}^{q} +
d_{1} \bar{d}_{2} s^{\bar{q}} \bar{s}^{-q} + \bar{d}_{1} d_{2}
s^{-\bar{q}} \bar{s}^{q} 
\vspace{3mm}\\
\ds \qquad + |d_{2}|^{2} s^{-\bar{q}} \bar{s}^{-q}
+|c_{1}|^{2} s^{r} \bar{s}^{\bar{r}} + \bar{c}_{1} c_{2} s^{-r} \bar{s}^{\bar{r}}
+c_{1} \bar{c}_{2} s^{r} \bar{s}^{-r} + |c_{2}|^{2} s^{-r} \bar{s}^{-\bar{r}}]
=0,
\vspace{3mm}\\
\ds \beta r (c_{1} s^{r} - c_{2} s^{-r}) +(d_{1} s^{\bar{q}} + d_{2} s^{-\bar{q}})
[ |d_{1}|^{2} \bar{s}^{q} s^{\bar{q}} + d_{1} \bar{d}_{2} s^{\bar{q}}
\bar{s}^{-q} + \bar{d}_{1} d_{2} s^{-\bar{q}} \bar{s}^{q} 
\vspace{3mm}\\
\ds \qquad +|d_{2}|^{2} s^{-\bar{q}} \bar{s}^{-q}
+|c_{1}|^{2} s^{r} \bar{s}^{\bar{r}}+\bar{c}_{1} c_{2} s^{-r} \bar{s}^{\bar{r}}
+ c_{1} \bar{c}_{2} s^{r} \bar{s}^{-\bar{r}} + |c_{2}|^{2} s^{-r} \bar{s}
^{-\bar{r}}] = 0,
\ea
\ee
and their respective complex conjugates. 
We require that system (7.14) is satisf\/ied for any value of $s$ and $\bar{s}$.
This means that the coef\/f\/icients of the successive powers of $s$ and $\bar{s}$
in equations (7.14) have to vanish.
As a result, we obtain a consistent system of algebraic 
equations for $c_{i}$ and $d_{i}$, $i=1,2$. 
The solutions of WE system (1.2) in this case take the form
\be
\ba{l}
\ds \psi_{1} =  D E \frac{\exp(3 \lambda (z + \bar{z})/2) \exp(\lambda(z - \bar{z})/2)}{E^{2} 
\exp(2 \lambda (z +\bar{z})) + 1},  
\vspace{3mm}\\ 
\ds \psi_{2} = D  \frac{\exp(\lambda (z + \bar{z})/2) \exp(\lambda(z- \bar{z})/2)}
{E^{2} \exp(2 \lambda (z + \bar{z})) + 1},
\ea
\ee
where $D = 2 (c+i)( \lambda E/(c^{2} + 1))^{1/2}$, and $c$, $\lambda$
and $E$ are arbitrary real constants of integration.
For $\lambda < 0$ solutions (7.15) are nonsingular and represent a
bump-type solutions. The associated surface is obtained from relations (1.3) as
\be
\ba{l}
\ds X_{1} = -2 \left(\left( c^{2} -1\right) \sin x + 2 c \cos x\right)
\left(1 + \frac{1}{2 \lambda}\right) \left(p^{2} + q^{2}\right),
\vspace{3mm}\\
\ds X_{2} =2 \left( 2 c \sin x - \left(c^{2} -1\right) \cos x\right)
\left(1 + \frac{1}{2 \lambda}\right) \left(p^{2}+ q^{2}\right),
\vspace{3mm}\\
\ds X_{3} = 4  \left(E^{2} e^{4 \lambda x} + 1\right)^{-1}.
\ea
\ee
Here, $x$ is the real part of $z$, and $p^{2}$ and $q^{2}$ are given by
\[
p^{2} = \frac{4 \lambda E^{3} e^{6 \lambda x}}{A (E^{2} e^{4 \lambda x}+1)^{2}},
\qquad
q^{2} = \frac{4 \lambda E e^{2 \lambda x}}{A( E^{2} e^{4 \lambda x} +1)^{2}},
\qquad A \in \mathbb R.
\]
One can eliminate quantities
$p$, $q$ and $x$ from (7.16) and obtain the 
following expression for the surface representing a catenoide
\[
X_{1}^{2} + X_{2}^{2} = 4 \left(\frac{\lambda}{A}\right)^{2} 
\left[(c-1)^{2} + 4 c^{2}\right] \left(1 + \frac{1}{2 \lambda}\right)^{2} (4 - X_{3}) X_{3}.
\]
The Gaussian curvature $K$ is constant and given by
\[
K = 4 \frac{\lambda^{2} E^{2}}{|D|^{2}}.
\]

\section{Final remarks} 

We have presented a variety of new approaches to the study
of the generalized WE system. They proved to be
particularly ef\/fective in delivering solutions from which it was
possible to derive explicit formulae for associated constant
mean curvature surfaces embedded in~$\mathbb R^{3}$.
One of the more interesting results of our analysis is the
observation that the WE system admits potential solutions.
This fact made it possible, for the f\/irst time, to deal
with the boundary value problem for the generalized WE system. We were
also able to construct potential multi-soliton solutions. 

It is worth noting that the treatments of the WE system proposed here
can be applied, with necessary modif\/ications, to more general 
cases of WE systems describing surfaces immersed in multi-dimensional
Euclidean and pseudo-Riemannian spaces. Such 
generalization of the WE system has been recently presented by Konopelchenko 
in [18] where, in particular, the explicit formulae for
minimal surfaces immersed in four-dimensional Euclidean
space $\mathbb R^{4}$ and $S^{4}$ have been derived.
An extention of our analysis to this case will be a 
subject of a future work.

\subsection*{Acknowledgements}
The authors would like to thank Professor B.~Konopelchenko (University of Lecce) for 
helpful discussions on this subject.
The research reported in this paper was partly
supported by the National Sciences and Engineering Research
Council of Canada and also by the Fonds FCAR du Gouvernment du Qu\'{e}bec.

\label{bracken-lp}

\end{document}